\RequirePackage{fix-cm}
\documentclass[oneolumn, twoside=off, hidelinks]{preprint}    

\smartqed 
\usepackage{graphicx}
\usepackage{hyperref}
\usepackage{tikz}
\usepackage{pgf}

\renewcommand{\vec}[1]{{\ensuremath{\boldsymbol{\mathrm #1}}}}
\newcommand{\ten}[1]{\ensuremath{\boldsymbol{\mathsf{#1}}}}
\newcommand{\tI}{{\ensuremath{\ten I}}}
\newcommand{\vdot}{\boldsymbol{\mathsf{\ensuremath\cdot}}}

\newcommand{\del}{\ensuremath{\nabla}}
\newcommand{\deld}{\ensuremath{\del\vdot}}
\newcommand{\lrp}[1]{\left( #1 \right)}

\newcommand{\BJ}{{\mathrm{BJ}}}
\newcommand{\FF}{{\mathrm{ff}}}
\newcommand{\PM}{{\mathrm{pm}}}
\newcommand{\RR}{\mathbb{R}}
\newcommand{\T}{\ensuremath{^{\mathrm T}}}

\usepackage{amsmath,amssymb}
\usepackage{color}
\usepackage{cite}

\usepackage[inner=3.5cm, outer=3.5cm,top=3cm,bottom=4cm]{geometry}

\begin{document}
\title{Validation and calibration of coupled porous-medium and free-flow 
problems using pore-scale resolved models}

\titlerunning{Validation and calibration of coupled flow models}        

\author{Iryna Rybak \and 
        Christoph Schwarzmeier \and \\
        Elissa Eggenweiler \and 
        Ulrich R\"ude 
}

\institute{
I. Rybak \and E. Eggenweiler \at
Institute of Applied Analysis and Numerical Simulation, \\
University of Stuttgart, \\ Pfaffenwaldring 57, 70569 
Stuttgart, Germany\\
\email{rybak@ians.uni-stuttgart.de, 
elissa.eggenweiler@mathematik.uni-stuttgart.de}
\and
C. Schwarzmeier \and U. R\"ude \at
Chair for System Simulation, \\
Friedrich-Alexander University Erlangen-N{\"u}rnberg, \\ Cauerstra{\ss}e 11, 91058 Erlangen, Germany \\
\email{christoph.schwarzmeier@fau.de, ulrich.ruede@fau.de}
\and U. R\"ude \at
CERFACS, 42 Avenue Gaspard Coriolis, 31057 Toulouse, France
}

\date{ }
\maketitle

\begin{abstract}
The correct choice of interface conditions and effective parameters for coupled 
macroscale free-flow and porous-medium models is crucial for a complete 
mathematical description of the problem under consideration and for 
accurate numerical simulation of applications. We consider single-fluid-phase systems described by the Stokes--Darcy model. Different sets of coupling 
conditions for this model are available. However, the choice of these conditions and effective model 
parameters is often arbitrary. We use large scale lattice Boltzmann 
simulations to validate coupling conditions by 
comparison of the macroscale simulations against pore-scale resolved 
models. We analyse two settings (lid driven cavity over a porous bed and 
infiltration problem) with different geometrical configurations 
(channelised and staggered distributions of solid grains) and different 
sets of interface conditions. Effective parameters for the macroscale 
models are computed numerically for each geometrical configuration. 
Numerical simulation results demonstrate the sensitivity of the coupled 
Stokes--Darcy problem to the location of the sharp fluid-porous interface, 
the effective model parameters and the interface conditions.

\end{abstract}

\section{Introduction}\label{intro}

Coupled flow systems containing a porous-medium domain and a free-flow 
region appear in various environmental and technical applications, such as 
surface-water/ groundwater flow, industrial filtration and water management
in fuel cells. The interaction between the flow regions is dominated by the 
interface driven processes and, due to their complexity, modelling such 
coupled flow systems is a challenging task.

Different spatial scales can be employed to investigate coupled free-flow 
and porous-medium systems. At the \emph{microscale} (pore scale) the pore 
structure is fully resolved (Fig.~\ref{fig:domain}, left) and the flow in the entire fluid domain 
(free-flow region and pore space in the porous medium) is described by the 
Navier--Stokes equations with the no-slip condition at the boundary of the 
solid inclusions. However, computing the microscale flow field is infeasible for practical 
applications as it
requires detailed information about porous-medium morphology and topography 
which is usually unknown. Even if information on the pore structure is 
available, e.g., obtained from tomographic analysis~\cite{Blunt_17, 
Wildenschild_Sheppard_13} or known in advance for artificially produced 
composite materials (GeoDict, \url{www.math2market.com}), performing 
microscale simulations is computationally highly demanding. 

From the \emph{macroscale} perspective, the whole system is described as 
two different continuum flow domains (free flow, porous medium) separated by an equi-dimensional transition region or 
a lower-dimensional (co-dimension one) sharp interface. Different mathematical models are 
usually used in the two flow regions. Thus, either coupling conditions at the sharp 
fluid-porous interface or equations that apply at the transition zone between 
the two flow systems are needed \cite{Beavers_Joseph_67, Goyeau_Lhuillier_etal_03, 
Nield_09, OchoaTapia_Whitaker_95, Saffman}. Correct specification of these 
conditions is essential for a complete and accurate mathematical 
description of flow and transport processes in compositional systems as 
well as for an accurate numerical simulation of applications 
\cite{Discacciati_Miglio_Quarteroni_02, Discacciati_Quarteroni_09, 
Jackson_Rybak_etal_12, Jaeger_Mikelic_09, Layton_Schieweck_Yotov_03}. The 
Navier--Stokes equations, which in case of low Reynolds numbers can be 
simplified to the Stokes equations, are typically applied to describe 
fluid flow in the free-flow domain, while Darcy's law or its extension is 
used for the flow through the porous medium. Depending on the flow 
regime, other models can be applied, e.g., shallow water equations or 
kinematic wave approximation of the Saint-Venant equation for surface 
flows~\cite{Dawson_08, 
Sochala_Ern_Piperno_09, Reuter_etal_19} and the Brinkman equation, the 
Forchheimer equation, the Richards equation or the full two-phase 
porous-medium flow models for subsurface flows~\cite{Brinkman_47, 
Cimolin_Discacciati_13, Mosthaf_Baber_etal_11, Rybak_etal_15}. Different 
sets of coupling conditions are used in the literature at the fluid-porous 
interface~\cite{Beavers_Joseph_67, OchoaTapia_Whitaker_95, Carraro_etal_15, 
Nield_09, Angot_etal_17}. Some authors distinguish between two 
qualitatively different 
flow directions: {\it near parallel flow}, where the velocity of the free 
fluid is much larger than the filtration velocity in the porous medium and 
{\it near normal flow}, where velocities in the two regions have 
comparable magnitudes. However, these conditions cannot be applied to a 
general flow situation. 

The macroscopic level of flow system description is the most common one  
used for numerical simulations of coupled systems in practical 
applications~\cite{Arbogast_Brunson_07, Cimolin_Discacciati_13, 
HanspalWaghodeNassehiWakeman, Iliev_Laptev_04, Mosthaf_Baber_etal_11, 
Riviere}. Microscale (pore-scale resolved) numerical simulations are 
typically used to validate the macroscale models. 

In this paper, we consider the Stokes--Darcy problem and analyse 
different coupling conditions at the fluid-porous interface. Most often, 
the conservation of mass, the balance of normal forces and the 
Beavers--Joseph--Saffman condition for the tangential velocity are used 
both for numerical analysis and simulations, e.g.,~\cite{Beavers_Joseph_67, 
Saffman, Goyeau_Lhuillier_etal_03, Nield_09, Discacciati_Quarteroni_09, 
GiraultRiviere, Layton_Schieweck_Yotov_03, Jaeger_Mikelic_09}. The 
Beavers--Joseph--Saffman condition that describes a tangential 
velocity slip at the fluid-porous interface was obtained for flows 
parallel to the porous layer. Nevertheless, this condition is often applied 
to other flow  
directions, e.g.,~\cite{HanspalWaghodeNassehiWakeman}. Although other interface 
conditions are available for the Stokes--Darcy 
problem, they are either valid only for specific flow 
situations~\cite{Carraro_etal_15} or contain parameters which are difficult to 
determine~\cite{ValdesParada_etal_09}. There were also 
several attempts to validate the coupled models and the interface 
conditions, e.g.,~\cite{Lacis_Bagheri_17, Zampogna_Bottaro_16, 
Bars_Worster_06, Yang_etal_19}.
However, to our knowledge no systematic study has been published.
In this paper, we consider different interface conditions for the 
coupled Stokes--Darcy problem. We validate and calibrate the macroscale 
coupled model using pore-scale resolved numerical simulations. These 
fully resolved numerical simulations are performed using the lattice Boltzmann 
method, a numerical method for simulating fluid dynamics that is especially 
suited for parallel computations with high spatial and temporal resolutions 
in complex geometries. It thus provides first-principle insights into the 
physical flow processes in both domains and the flow behaviour at the 
interface region. 

The paper is organised as follows. In Section~\ref{sec:geom} we provide 
the geometrical setting and the flow system description. The 
macroscale coupled model with different sets of interface conditions is 
presented in Section~\ref{sec:model}. The calculation of effective 
parameters for the macroscale model by means of homogenisation is 
discussed in Section~\ref{sec:eff-param} and the lattice Boltzmann method 
is described in Section~\ref{sec:LBM}. Numerical simulation results for 
model validation and calibration are given in 
Section~\ref{sec:validation}. Finally, the conclusion is presented in 
Section~\ref{sec:conclusion}.

\section{Geometrical setting}\label{sec:geom}

\begin{figure}[!h]
\centerline{\includegraphics[width=0.4\textwidth]{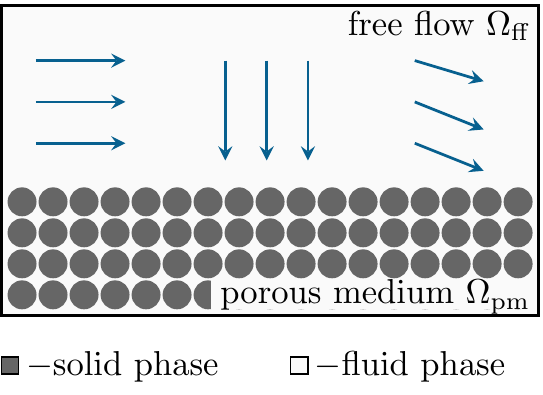}\;
\includegraphics[width=0.4\textwidth]{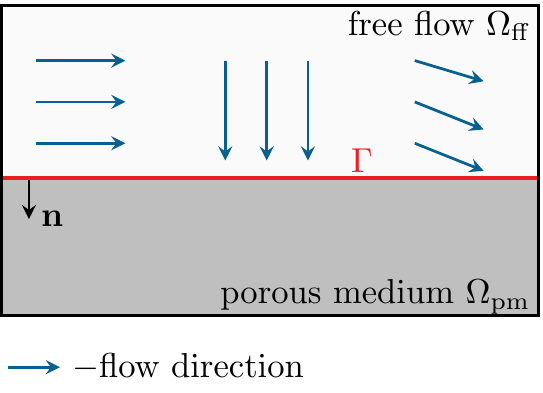}}
\caption{Schematic representation of the flow system at the microscale 
(left) and the macroscale (right).}
\label{fig:domain}
\end{figure}

We consider a free-flow region $\Omega_\FF \subset \RR^2$ and a porous 
medium $\Omega_{\text{pm}} \subset \RR^2$ with periodic arrangement of the 
solid obstacles (Fig.~\ref{fig:domain}, left). The fluid is considered to be present in a single phase and to fully saturate the porous medium. Additionally, it is assumed to be incompressible and to have a constant viscosity. The flow is considered as being laminar and supposing a single chemical species, no compositional effects need to be modelled. The temperature is assumed constant, such that it is not required to model any energy balance equation.

At the microscale, we solve the Stokes equations in the free-flow domain 
and in the pore space of the porous medium, where no-slip boundary 
conditions are applied on the solid-fluid boundary. At the macroscale, the 
system is described as two different continuum flow domains separated by 
the sharp interface $\Gamma$ (Fig.~\ref{fig:domain}, right).
We solve the Stokes and Darcy equations in the two flow domains and impose 
different interface conditions on the fluid-porous interface 
(Sect.~\ref{sec:model}). These interface conditions are validated for 
different porous-medium morphologies and different flow directions. 

\section{Macroscale coupled model formulation}\label{sec:model}

Fluid flow in the free-flow region is described by the Stokes 
equations, the flow through the porous medium is described by Darcy's 
law. Different sets of interface conditions are imposed on the sharp 
fluid-porous interface. Effective parameters such as permeability, 
porosity and boundary layer constants are needed to describe the 
porous-medium system and the interface conditions. These parameters are 
calculated in 
Section~\ref{sec:validation} for the geometrical configurations 
studied in the paper.

\subsection{Free-flow model}

The fluid is considered to be incompressible and therefore the mass
conservation equation becomes 
\begin{equation}\label{eq:1p1c-FF-mass}
\deld \vec v = 0 \qquad \text{in} \;\; \Omega_\FF.
\end{equation}
The Stokes equation describes the fluid flow 
\begin{equation}\label{eq:1p1c-FF-momentum}
- \deld \ten T(\vec v,p) - \rho \vec g = 0
\qquad \text{in} \;\; \Omega_\FF,
\end{equation}
where $\rho$ is the fluid density, $\vec v$ is the fluid velocity,
$p$ is the pressure, $\vec g$ is the gravitational acceleration,
$\ten T(\vec v,p) = 2 \mu \ten D\lrp{\vec v} - p \tI$ is the stress
tensor, $\mu$ is the dynamic viscosity, $\displaystyle \ten D\lrp{\vec v} =
\tfrac12 \lrp{\del \vec v + \lrp{\del \vec v}^{\T}}$ is the rate of
strain tensor and $\tI$ is the identity tensor. 

On the external boundary of the free-flow domain $\Gamma_\FF = \partial 
\Omega_\FF \setminus\Gamma$, the following boundary conditions 
are imposed
\begin{equation}\label{eq:FF-BC}
\begin{split}
\vec v &= \overline {\vec v}_\FF \hspace*{3.5mm} \text{on} \;\; \Gamma_{D, \FF}, 
\\
\ten T(\vec v, \, p) \vdot \vec n_\FF &= \overline {\vec 
h}_\FF  \hspace*{3.5mm} \text{on} \;\; \Gamma_{N, \FF},
\\
\partial \vec v_\FF/\partial \vec n_\FF &= 0  \hspace*{5.25mm} \; \text{on} \;\; 
\Gamma_{\text{out}, \FF},
\end{split}
\end{equation}
where $\vec n_\FF$ is the unit outward normal vector from domain 
$\Omega_\FF$ on its boundary and $\Gamma_\FF = \Gamma_{D, \FF} \cup 
\Gamma_{N, \FF} \cup \Gamma_{\text{out}, \FF}$.

\subsection{Porous-medium model}

We consider Darcy's flow equations in the porous-medium domain 
\begin{equation}\label{eq:PM-mass}
\deld \vec v = q \qquad \text{in} \;\; 
\Omega_\PM, 
\end{equation}
where Darcy's velocity is given by
\begin{equation}\label{eq:PM-Darcy}
\vec v = -\frac{\ten K}{\mu} \lrp{\del p - \rho \vec g}.
\end{equation}
Here $\vec v$ is the fluid velocity through the porous medium, $p$ 
is the fluid pressure, $\ten K$ is the intrinsic permeability of the 
porous medium, $\mu$ is the dynamic viscosity of the fluid, $\rho$ is the 
fluid density, $q$ is the source term and $\vec g$ is the gravitational 
acceleration. The permeability tensor $\ten K$ is symmetric, 
positive definite and bounded. Its values are computed 
numerically for the considered settings in 
Section~\ref{sec:validation}. 
In this paper, we neglect gravitational effects setting $\vec g=\vec 0$.

We prescribe the following boundary conditions on the external boundary of 
the porous-medium domain $\Gamma_\PM = \partial 
\Omega_\PM \setminus \Gamma$: 
\begin{equation}\label{eq:PM-BC}
\begin{split}
p &= \overline p_\PM \quad \text{on} \;\; \Gamma_{D, \PM}, \\
\vec v \vdot \vec n_\PM &= \overline v_\PM \quad \text{on} \;\; \Gamma_{N, 
\PM},
\end{split}
\end{equation}
where $\vec n_\PM$ is the unit outward 
normal vector from domain $\Omega_\PM$ on its boundary, $\Gamma_\PM = 
\Gamma_{D, \PM} \cup \Gamma_{N, \PM}$, $\Gamma_{D, \PM} \cap \Gamma_{N, 
\PM}=\emptyset$, and $\Gamma_{D, \PM} \neq \emptyset$.

\subsection{Interface conditions}

In addition to the boundary conditions prescribed on the external boundary 
of the coupled domain, a set of interface conditions has to be defined on 
$\Gamma$ in order to obtain a closed problem formulation. Different sets 
of interface conditions have been proposed in the literature depending on the flow 
direction. 

\subsubsection{Flow parallel to the interface}

When the free flow is almost parallel to the porous medium, the 
conservation of mass, the balance of normal forces and a variation of the 
Beavers--Joseph condition is typically applied.

The {\it conservation of mass} across the interface $\Gamma$ is given by 
\begin{equation}\label{eq:IC-mass}
\vec v_\FF \vdot \vec n = \vec v_\PM \vdot \vec n \qquad \text{on} 
\;\Gamma,
\end{equation}
where $\vec n$ is the unit normal vector at the fluid-porous interface 
$\Gamma$ pointing outward from the free-flow domain $\Omega_\FF$ 
(Fig.~\ref{fig:domain}, right).

The {\it balance of normal forces} at the interface $\Gamma$ is 
\begin{equation}\label{eq:IC-momentum}
-\vec n_\FF \vdot \ten T\lrp{\vec v_\FF, p_\FF} \vdot \vec n_\FF = p_\PM 
\qquad \text{on} \; \Gamma,
\end{equation}
where $\vec n_\FF = \vec n = -\vec n_\PM$. 

There exist several possibilities for the interface condition on the 
tangential component of the velocity~\cite{Nield_09}. The {\it 
Beavers--Joseph} interface condition \cite{Beavers_Joseph_67} for the 
difference between the tangential component of the free-flow velocity and 
the porous-medium velocity can be written as
\begin{equation}\label{eq:IC-BJ}
\lrp{\vec v_\FF -\vec v_\PM} \vdot \vec \tau_j + \frac{\sqrt{\ten 
K}}{\alpha_{\text{BJ}}} \lrp{ \del \vec v_\FF \vdot \vec n_\FF} \vdot \vec
\tau_j = 0 \quad \text{on} \; \Gamma,
\end{equation}
where $\alpha_{\text{BJ}}>0$ is the Beavers--Joseph parameter, $\vec 
\tau_j$ is a unit vector tangential to the interface, $j = 1,\ldots,d-1$, 
and $d$ is the number of space dimensions.

Saffman~\cite{Saffman} simplified the Beavers--Joseph condition 
\eqref{eq:IC-BJ} as follows
\begin{equation}\label{eq:IC-BJS}
\vec v_\FF \vdot \vec \tau_j + \frac{\sqrt{\ten K}}{\alpha_\BJ}
\lrp{ \del \vec v_\FF \vdot \vec n_\FF} \vdot \vec \tau_j
= 0 \quad \text{on} \; \Gamma,
\end{equation}
arguing that the tangential velocity of a fluid through a porous medium is 
small and thus can be neglected.

Jones~\cite{Jones_73} considered the symmetrised rate of
strain tensor $\ten D$ in the 
Beavers--Joseph--Saffman interface condition 
\begin{equation}\label{eq:IC-BJJ}
\vec v_\FF \vdot \vec \tau_j + \frac{\sqrt{\ten 
K}}{\alpha_{\text{BJ}}} \vec n_\FF \vdot \ten D\lrp{\vec v_\FF} \vdot \vec
\tau_j = 0 \qquad \text{on} \; \Gamma.
\end{equation}

Besides the classical interface conditions \eqref{eq:IC-mass}, 
\eqref{eq:IC-momentum} and \eqref{eq:IC-BJ}--\eqref{eq:IC-BJJ}, which have 
been formulated based on heuristic arguments and were justified 
experimentally, in \cite{Jaeger_etal_01} the following interface 
conditions on $\Gamma$ were proposed
\begin{equation}\label{eq:IC-Mik-parallel}
\begin{split}
v_{\FF,2} &= 0, \\
p_{\text{ff}} + \mu 
C_{\omega}^{\text{bl}} \frac{\partial v_{\text{ff}, 1}}{\partial 
x_2} &= p_{\text{pm}},\\
v_{\text{ff}, 1} - \varepsilon C_1^{\text{bl}} \frac{\partial 
v_{\text{ff},1}}{\partial x_2} &= 0,
\end{split} 
\end{equation}
where $\vec v = (v_1,v_2)$ and $\varepsilon$ is the ratio of the length 
scales (Sect.~\ref{sec:eff-param}). These coupling conditions are derived 
mathematically using homogenisation theory and boundary layer theory for 
flows parallel to the fluid-porous interface $\Gamma$. 
The first condition in \eqref{eq:IC-Mik-parallel} differs from the 
classical conservation of mass condition \eqref{eq:IC-mass}. However, it 
seems to be reasonable for flows parallel to the porous layer, because in 
this case the normal component of the Darcy velocity is of order 
$\varepsilon^2$ (see Eq.~\eqref{eq:asymexp}) and thus can be neglected. 
The second condition in 
\eqref{eq:IC-Mik-parallel} links pressure fields in the porous domain and 
the free-flow domain, and is a variation of the classical balance of normal 
forces~\eqref{eq:IC-momentum}. The last equation in 
\eqref{eq:IC-Mik-parallel} is a modification of the 
Beavers--Joseph--Saffman law~\cite{Beavers_Joseph_67, Saffman} with 
$\displaystyle \alpha_\BJ^{-1} \ten{K}^{1 \slash 2} = - \varepsilon 
C_1^{\text{bl}}\tI$. The effective coefficients $\displaystyle 
 C_1^{\text{bl}}$ and $\displaystyle C_\omega^{\text{bl}}$ are determined 
through an auxiliary boundary layer problem. These coefficients and the 
parameter $\varepsilon$ are provided for different geometrical 
configurations in Section~\ref{sec:validation}.

\subsubsection{Flow perpendicular to the interface} 

For flows almost perpendicular to the porous medium, the following 
effective interface laws on $\Gamma$ are derived in \cite{Carraro_etal_15} 
using the theory of homogenisation and boundary layers
\begin{equation}
\begin{split}
v_ {\FF, 2} &=  v_{\PM, 2}, 
\\
v_ {\FF, 1} &= M_1^{\text{bl}} \frac{\partial 
p_\PM}{\partial x_2},
\\
p_\PM &=0,
\end{split}
\label{eq:IC-Mik-normal}
\end{equation}
where $M_1^{\text{bl}}$ is a boundary layer corrector constant. However, 
these conditions are valid only for a very specific setting (given 
geometrical configuration and boundary conditions) and cannot be applied to 
a general infiltration problem. There is a lack of physically consistent 
interface conditions for flows perpendicular to the porous layer and 
moreover for arbitrary flow directions.

\section{Effective properties of the coupled system}\label{sec:eff-param}

\begin{figure}[ht!]
\centering
\includegraphics[scale=0.9]{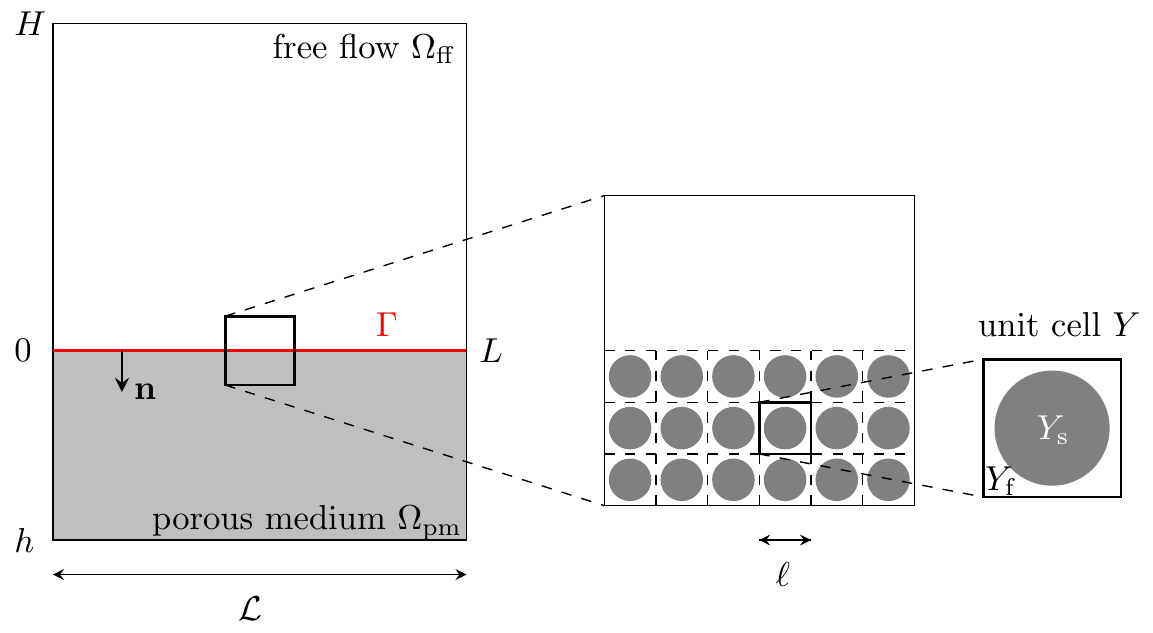}
\caption{Separation of scales: $\ell$ - characteristic pore size, 
$\mathcal{L}$ - macroscopic domain length.}
\label{fig:unitcell}
\end{figure}

To obtain effective properties for coupled macroscale models 
(effective permeability and boundary layer constants), we use the 
theory of homogenisation and boundary layers \cite{Carraro_etal_15, 
Hornung_97, Jaeger_Mikelic_09}. Therefore, we consider periodic porous 
media and assume the separation of length scales $\varepsilon = 
\ell/\mathcal{L} \ll 1$, where $\ell$ is the characteristic pore size and 
$\mathcal{L}$ is the length of the domain of interest 
(Fig.~\ref{fig:unitcell}).

We consider a coupled domain $\Omega = \Omega_\FF \cup \Omega_\PM$, 
consisting of a free-flow region $\Omega_\FF \subset \RR^2$ and a porous 
medium $\Omega_{\text{pm}} \subset \RR^2$ with  periodic arrangement of 
the solid obstacles (Fig.~\ref{fig:domain}). The flow region 
$\Omega^{\varepsilon} = \Omega_{\text{ff}} \cup \Gamma 
\cup \Omega_{\text{pm}}^\varepsilon$ consists of the free-flow 
domain $\Omega_\FF$, the permeable interface $\Gamma$ and the pore 
space of the porous-medium domain $\Omega_{\text{pm}}^\varepsilon$. 
The porous-medium geometry is defined by a periodic repetition of the scaled unit 
cell $Y=(0,1)^2$ in $\Omega_\PM$ (Fig.~\ref{fig:unitcell}).

\subsection{Microscale model}\label{sec:pore-scale}

The fluid flow in $\Omega^{\varepsilon}$ is governed by the non-dimensional 
steady Stokes equations
\begin{equation}\label{eq:microscale}
\begin{split}
- \mu \Delta \textbf{v}^{\varepsilon}+ \nabla p^{\varepsilon} &= 0 \qquad 
\textnormal{in } \Omega^{\varepsilon},
\\
\deld \textbf{v}^{\varepsilon} &= 0 \qquad \textnormal{in } 
\Omega^{\varepsilon},
\end{split}
\end{equation}
completed with the no-slip condition on the fluid-solid interface 
\begin{equation}\label{eq:no-slip}
\vec v^\varepsilon = \vec 0 \quad \text{on } \partial \Omega^\varepsilon 
\setminus \partial \Omega,
\end{equation}
and the appropriate boundary conditions on the external 
boundary $\partial \Omega$. 

For the lid driven cavity over the porous bed described in 
Section~\ref{sec:lid-driven}, we consider 
\begin{equation}\label{eq:microscale-BC}
\begin{split}
\textbf{v}^{\varepsilon} &= (0,0) \qquad \textnormal{on } \partial 
\Omega \setminus\{x_2 = H\},
\\
\textbf{v}^{\varepsilon} &= (1,0) \qquad \textnormal{on } \{x_2 = H\},
\end{split}
\end{equation}
and for the infiltration problem validated in 
Section~\ref{sec:infiltration}, we have
\begin{align}\label{eq:microscale-BC2}
\begin{split}
&\vec v^{\varepsilon} = (0,0) \hspace*{21.75mm}  \textnormal{on } \{x_1 
= 0\} \cup \{x_1 = L\}, \\
&\textbf{v}^{\varepsilon} = (0,-\sin(\pi x_1)) \hspace*{7.5mm} \textnormal{on 
} \{x_2 = H\}, \\
& \left( \nabla \textbf{v}^{\varepsilon} - p^{\varepsilon} \textbf{I} 
\right) \vec n = -p_b\vec n  \hspace*{3.5mm}\textnormal{on } 
\{x_2 = h\},
\end{split}
\end{align}
where $p_b = 100$.

\subsection{Homogenisation}\label{sec:homogenization}
We obtain the corresponding macroscale formulation of the coupled system 
by studying the behaviour of the solutions to the microscopic 
problem~\eqref{eq:microscale}--\eqref{eq:microscale-BC} or 
\eqref{eq:microscale}, \eqref{eq:no-slip} and \eqref{eq:microscale-BC2} 
when $\varepsilon 
\rightarrow 0$. In the limit, the equations in the free-flow region 
$\Omega_\FF$ remain unchanged and we need to homogenise the porous 
region. 
To derive the limit problem in $\Omega^{\varepsilon}_{\text{pm}}$, we use 
the idea of homogenisation with two-scale asymptotic 
expansions~\cite[chap.\,1 and 3]{Hornung_97}:  
\begin{equation}\label{eq:asymexp}
\begin{split}
\textbf{v}^{\varepsilon} (\textbf{x}) &= \varepsilon^2 \textbf{v}_0 (\textbf{x}, \textbf{y}) + \varepsilon^3 \textbf{v}_1 (\textbf{x}, \textbf{y}) + \ldots 
\\
p^{\varepsilon} (\textbf{x}) &= p_0 (\textbf{x}, \textbf{y}) + \varepsilon p_1 (\textbf{x}, \textbf{y}) + \ldots 
\end{split}
\end{equation}
where $\vec y = \vec x \slash \varepsilon$ is the fast variable. We follow 
the classical procedure of homogenisation to derive Darcy's law 
\eqref{eq:PM-Darcy} as the upscaled problem in $\Omega_{\text{pm}}$
completed by the appropriate boundary conditions. 

The permeability tensor $\ten{\widetilde{K}}$ is given by 
\begin{equation}\label{eq:permeability}
\begin{split}
\ten{\widetilde{K}} = (\widetilde k_{ij})_{i,j = 1,2} = \int_{Y_\text{f}} 
\text{w}_{i}^j \ \text{d} 
\textbf{y},  
\end{split}
\end{equation}
where $\textbf{w}^{j} = (w_1^j, w_2^j)$ and $\pi^j$ are the solutions to 
the following cell 
problems for $j=1,2$:  
\begin{equation}\label{eq:cell-prob}
\begin{split}
- \Delta_ {\textbf{y}} \textbf{w}^{j}  + \nabla_{\textbf{y}} \pi^j &= 
\textbf{e}_j \quad \text{in $Y_\text{f}$},  
\\
\del_ {\textbf{y}} \vdot \textbf{w}^{j}  &= 0 \quad \;  \; \text{in 
$Y_\text{f}$}, 
\\
\textbf{w}^{j} &= 0 \quad \;  \; \text{on $\partial Y_\text{f} \setminus \partial Y$}, 
\\
\{ \textbf{w}^{j}, \pi^j \}& \text{ is 1-periodic in } \textbf{y}.
\end{split}
\end{equation}
The unit cell $Y$ is presented in Figure~\ref{fig:unitcell}.
To obtain the effective permeability $\ten{K}$, the tensor 
$\ten{\widetilde K}$ has to be scaled with $\varepsilon^2$.

\section{Lattice Boltzmann method}\label{sec:LBM}
The lattice Boltzmann method (LBM) is a modern approach for 
simulating fluid dynamics. While the LBM historically evolved from lattice 
gas automata \cite{Wolf-Gladrow_00}, it can be regarded as a special  
discretisation of the continuous Boltzmann equation~\cite{He_97} that 
describes the dynamics of a gas on a mesoscopic scale. On a macroscopic 
scale, however, the Boltzmann equation leads to the equations of fluid 
dynamics, i.e., mass, momentum and energy conservation. Therefore, by 
solving the Boltzmann equation for certain scenarios, one can also obtain 
the corresponding solutions to the Navier--Stokes 
equations~\cite{Krueger_LBM}.

\subsection{Advantages}\label{subsec:lbm-adv}
In comparison to conventional methods for simulating fluid dynamics, such 
as the finite difference, finite volume or finite element method, the 
LBM's 
main advantages lay in its ease of implementation and its suitability for 
complex flow scenarios. As the most compute-intensive operations are 
local, 
i.e., restricted to neighbouring lattice cells, the LBM is inherently 
applicable for massively parallel computing. Furthermore, the LBM is well-suited for simulating particulate 
flows~\cite{Rettinger_17, Kuron_16}, multi-phase flows~\cite{Huang_15} and flows through 
complicated geometries such as porous media~\cite{Fattahi_etal_16, 
Fattahi_etal_16b}.

\subsection{Description}\label{subsec:lbm-descr}
Our introduction to the LBM follows~\cite{Krueger_LBM}; the interested 
reader is referred to this reference and the references therein. 
The LBM's most fundamental quantities are the so-called particle distribution 
functions $f_{i}(\vec{x},t)$ that represent the density of virtual particles with 
velocity $\vec{c}_{i}=(c_{ix}, c_{iy})$ at position $\vec{x}$ and 
time $t$. The index $i \in \left[0; q\right)$ in $f_{i}$ refers to the 
corresponding index of velocity $\vec{c}_{i}$ in a discrete set of 
velocities $\{\vec{c}_{i}, w_{i}\}$ with weighting coefficients $w_{i}$.

Velocity sets are denoted as D$d$Q$q$ with $d$ being the number 
of spatial dimensions and $q$ referring to the velocity set's number of velocities. 
Velocity sets are a trade-off between accuracy and memory consumption. 
Different velocity sets and their corresponding coefficients can be found 
in the LBM literature such as~\cite{Krueger_LBM}. Here we limit ourselves to two spatial dimensions such that we use the 
well-established D2Q9 velocity set. In each velocity set, the constant 
$c_{s}$ determines the relation between the fluid pressure $p$ and the fluid density $\rho$ as 
$p = c_{s}^{2}\rho$. Due to this relation, $c_{s}$ is called 
the lattice speed of sound. In most velocity sets used to solve the 
Navier--Stokes equations, it can be expressed as $c_s = 
\sqrt{1/3} \Delta x/\Delta t$.

The distribution function $f_{i}$ is only defined at discrete points 
$\vec{x}$ in a square lattice with lattice spacing $\Delta x$. Since the LBM's origins lay in lattice gas automata, 
these discrete points are also often referred to as lattice cells. Similarly, 
$f_{i}$ is only defined at discrete times $t$ with time step $\Delta t$.
As it is common in the LBM literature, we use lattice units which is an 
artificial set of units that is scaled such that the relations 
$\Delta x=1$ and $\Delta t=1$ hold.

In kinetic theory, one obtains the macroscopic mass and momentum density 
by taking the moments of the Boltzmann equation's distribution functions. 
Similarly, one can find the mass density $\rho$ and 
the momentum density $\rho\vec{v}$ with fluid velocity $\vec{v}$ by
\begin{equation}\label{eq:lbm-rho-u}
\rho(\vec{x},t) = \sum_{i} f_{i}(\vec{x},t),
\quad
\rho\vec{v}(\vec{x},t) = \sum_{i} \vec{c}_{i} f_{i}(\vec{x},t).
\end{equation}

\begin{figure}[!h]
	\centerline{\includegraphics[width=0.35\textwidth]{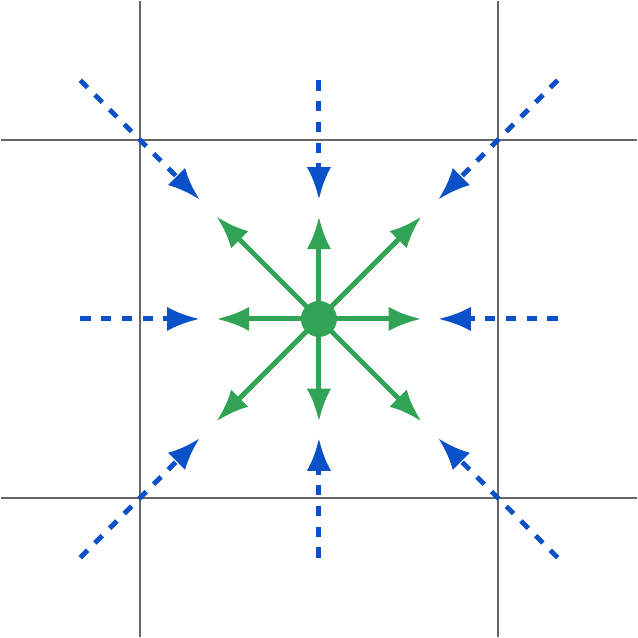}\;
		\hspace{0.035\textwidth}\includegraphics[width=0.35\textwidth]{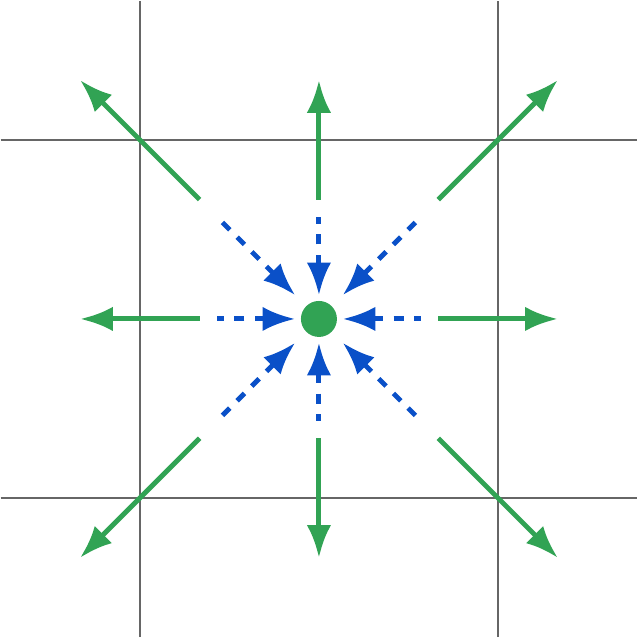}}
	\caption{Schematic representation of the lattice Boltzmann streaming and collision step for the D2Q9 velocity set. After particles have collided, they are streamed to neighbouring lattice cells (left). After streaming to neighbouring lattice cells, another particle collision occurs (right).}
	\label{fig:lbm-streamcollide}
\end{figure}

The lattice Boltzmann equation
\begin{equation}\label{eq:lbm-eq}
f_{i}(\vec{x} + \vec{c}_{i} \Delta t, t + \Delta t) = 
f_{i}(\vec{x}, t) + 
\Theta_{i}(\vec{x}, t)
\end{equation}
describes the advection of particles $f_{i}(\vec{x},t)$ that travel with 
velocity $\vec{c}_{i}$ to a neighbouring lattice cell 
$\vec{x}+\vec{c}_{i}\Delta t$ in the next time step $t+\Delta t$ (Fig.~\ref{fig:lbm-streamcollide}, left).
After having moved from one lattice cell to another, the particles are affected by a 
collision operator $\Theta_{i}(\vec{x},t)$ that models particle collisions by 
redistributing particles among $f_{i}$ (Fig.~\ref{fig:lbm-streamcollide}, right).

The simplest collision operator that can be used for solving the 
Navier--Stokes equations is the Bhatnagar--Gross--Krook (BGK) operator. It linearly relaxes particle distributions towards an equilibrium $f_{i}^{\text{eq}}$ by
\begin{equation}\label{eq:lbm-bgk}
\Theta_{i}(\vec{x},t) = - \frac{f_{i}(\vec{x},t)-f_{i}(\vec{x},t)^{\text{eq}}}{\tau}\Delta t,
\end{equation}
with a relaxation rate determined by the relaxation time $\tau$. In contrast to other collision operators, the BGK 
operator is based on only one relaxation time and it is therefore also called single-relaxation-time (SRT) 
collision operator. Due to its simplicity, the BGK collision operator is often the standard choice for LBM simulations. 
However, in~\cite{Pan_06}, it was found to be subject to a non-physical 
dependence between the fluid viscosity and the boundary 
locations, which is especially problematic in scenarios with 
large boundary areas such as flow through porous media.
Therefore, for the LBM simulations in this paper, we use the more accurate two-relaxation-time (TRT) collision 
operator~\cite{Ginzburg_08a,Ginzburg_08b}.

The particle distribution functions' equilibria
\begin{equation}\label{eq:lbm-feq}
f_{i}^{\text{eq}}(\vec{x},t) = 
w_{i}\rho\left(
1+
\frac{\vec{v}\cdot\vec{c}_{i}}{c_{\text{s}}^{2}}+
\frac{(\vec{v}\cdot\vec{c}_{i})^{2}}{2c_{\text{s}}^{4}}-
\frac{\vec{v}\cdot\vec{v}}{2c_{\text{s}}^{2}}
\right)
\end{equation}
are derived by a Hermite polynomial expansion of the continuous Maxwell-Boltzmann distribution~\cite{Krueger_LBM}.

The moments of $f_{i}^{\text{eq}}$ are identical to those of the 
corresponding $f_{i}$, i.e., 
$\Sigma_{i}f_{i}^{\text{eq}}=\Sigma_{i}f_{i}=\rho$ and 
$\Sigma_{i}\vec{c}_{i}f_{i}^{\text{eq}}=\Sigma_{i}\vec{c}_{i}f_{i}=\rho\vec
{v}$.  
Since $f_{i}^{\text{eq}}$ only depends on the local quantities density $\rho$ and fluid velocity 
$\vec{v}$ (Eq.~\eqref{eq:lbm-rho-u}), global information exchange is not needed and the LBM is therefore well-suited 
for parallel computing.

Using the Chapman--Enskog analysis~\cite{Krueger_LBM}, one can show that the lattice 
Boltzmann equation~\eqref{eq:lbm-eq} has a macroscopic behaviour that converges asymptotically to the solutions of the Navier--Stokes equations.
The fluid's kinematic shear viscosity $\nu$ is then a function of the 
relaxation time $\tau$ and obtained by
\begin{equation}\label{eq:lbm-viscosity}
\nu = c_{\text{s}}^{2}\left(\tau-\frac{\Delta t}{2}\right).
\end{equation}

\section{Model validation and calibration}\label{sec:validation}

In this section, we validate and calibrate coupled Stokes--Darcy 
models using pore-scale resolved simulations for different geometrical 
configurations and different flow regimes.

We consider the free-flow region $\Omega_\FF = [0,1] \times [0,1]$ 
and the porous medium $\Omega_\PM = [0,1] \times [-0.5,0]$ and study two 
test cases (lid driven cavity over a porous bed and infiltration problem). 
The porous medium is set to be periodic for both test cases such 
that we can compute the effective parameters easily using homogenisation 
theory. Moreover, we analyse two different periodic geometries. 

The porosity of the medium is taken to be $\phi = 0.4$. Therefore, 
the radius of solid inclusions is given by $r = \varepsilon \sqrt{(1 - 
\phi)/\pi}$. We solve the unit cell problems \eqref{eq:cell-prob} and 
compute the effective permeability tensor for different geometrical 
settings according to equation \eqref{eq:permeability}. To compute 
the boundary layer constants, we solve the Stokes problem in a perforated 
vertical stripe containing unit cells (see e.g. \cite{Carraro_etal_15}). 
For numerical simulations we use the software package \textsc{FreeFEM++}  
with P2/P1 finite elements~\cite{MR3043640}. The dynamic viscosity of the 
fluid is scaled for both test cases to $\mu = 1$. 

The macroscale problem is discretised using the finite volume method on 
staggered grids~\cite{VersteegMalalasekra}. The computational  
domains $\Omega_\FF$ and $\Omega_\PM$ are partitioned into equal blocks 
with grid size $h=1/800$, and the grids are conforming at the interface.

We apply large scale simulations using the LBM described in 
Section~\ref{sec:LBM} to compute the microscale (pore-scale) velocity 
field 
and compare the macroscale simulations with different sets of interface 
conditions against these pore-scale resolved models. 

The LBM simulations are implemented and performed with the 
open-source software-framework 
\textsc{waLBerla} \cite{Godenschwager_2013, Schornbaum_18} (\url{www.walberla.net}).
We use the $D2Q9$ velocity set and the TRT collision operator. We set the TRT collision operator's even order 
relaxation time $\tau = 1$ and choose its \textit{magic parameter} $\Lambda=3/16$~\cite{Ginzburg_08a,Ginzburg_08b}. 
With equation \eqref{eq:lbm-viscosity}, the fluid's kinematic viscosity in lattice units becomes $\nu=1/6$. We consider 
the 
simulation to have reached steady state when the relative temporal change in the absolute value of the fluid's maximum 
velocity in $x_2$-direction is below $10^{-4}$.

In the following sections, we analyse the sensitivity of the coupled 
Stokes--Darcy models to the location of the fluid-porous interface 
$\Gamma$, the choice of the 
effective parameters and the interface conditions.

\subsection{Lid driven cavity over the porous bed}\label{sec:lid-driven}

For the lid driven cavity problem, we consider two geometrical configurations, a channelised one schematically presented in Figure~\ref{fig:lid-driven} and a staggered arrangement of the solid inclusions presented in Figure~\ref{fig:infiltration}. In both cases there are 40 
solid inclusions in $x_1$-direction and 20 inclusions in 
$x_2$-direction in the porous domain $\Omega_\PM$.
We define the Reynolds number
\begin{equation}
Re=\frac{v_{\text{lid}} \, \mathcal{L}}{\nu}
\end{equation}
using the $x_1$-component of the velocity at the lid $v_{\text{lid}}$, the length of the domain $\mathcal{L}$ in 
$x_1$-direction (Fig.~\ref{fig:unitcell}) and the kinematic fluid 
viscosity $\nu$. 

\subsubsection{Channelised geometrical configuration}

\begin{figure}[!h]
\centering
\includegraphics[width=0.5\textwidth]{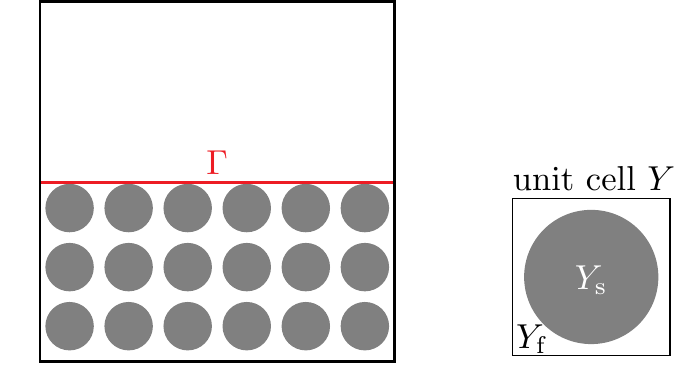}
\caption{Schematic illustration of the microscale setting and the 
location of the sharp fluid-porous interface $\Gamma$ (left) and the 
unit cell (right) for the channelised lid driven cavity problem.}
\label{fig:lid-driven}
\end{figure}

In this section, we consider the channelised geometrical configuration 
(Fig.~\ref{fig:lid-driven}). For this geometry, we solve the cell problems 
 \eqref{eq:cell-prob}, compute 
the permeability tensor according to \eqref{eq:permeability} and 
obtain 
\begin{equation}\label{eq:k1}
\ten{\widetilde{K}} =
\begin{pmatrix}
5.671\cdot 10^{-4} & 0 \\
0 & 5.671\cdot 10^{-4}
\end{pmatrix}.
\end{equation}
The effective permeability tensor of the porous medium is given by $\ten K 
= \varepsilon^2 \ten{\widetilde{K}}$, where the actual scale ratio is $\varepsilon = 1/40$.

The corresponding macroscale model \eqref{eq:1p1c-FF-mass}, 
\eqref{eq:1p1c-FF-momentum}, \eqref{eq:PM-mass}, \eqref{eq:PM-Darcy}  is 
closed by the following boundary conditions on the external boundary 
\begin{equation}\label{eq:lid-bc}
\begin{split}
&\vec v_{\FF} = (1,0) \hspace*{9.5mm}  \textnormal{on } \{x_2 
= 1\},
\\
&\vec v_{\FF} = (0,0) \hspace*{9.5mm}  \textnormal{on } (\{0 \} \cup \{1\} ) 
\times (0,1),
\\
&\vec v_{\PM} \vdot \vec n_\PM = 0 \hspace*{6.2mm}  \textnormal{on } \{x_2 
= -0.5\} 
\\
&\vec v_{\PM} \vdot \vec n_\PM = 0 \hspace*{6.2mm} \textnormal{on } (\{0 \} 
\cup 
\{1\} ) \times (-0.5,0).
\end{split}
\end{equation}
\vspace*{3mm}

In the LBM simulation, the whole domain is resolved by $2400 \times 3600$ 
lattice cells such that a solid inclusion has 
a diameter of about 52 lattice cells. We apply simple bounce back boundary conditions~\cite{Krueger_LBM} at the solid 
inclusions and at the external boundaries of the domain with respect 
to~\eqref{eq:lid-bc}. In order to approximate Stokes flow, 
we set $Re=0.01$.

Two sets of interface conditions are analysed for the macroscale problem: 
(i) the classical conditions \eqref{eq:IC-mass}, \eqref{eq:IC-momentum}, 
\eqref{eq:IC-BJJ} and (ii) the interface conditions derived using 
homogenisation \eqref{eq:IC-Mik-parallel}.

\subsubsection*{(i) Classical interface conditions}

In this section, we study the sensitivity of the coupled macroscale 
Stokes--Darcy (SD) model \eqref{eq:1p1c-FF-mass}, 
\eqref{eq:1p1c-FF-momentum}, \eqref{eq:PM-mass}, \eqref{eq:PM-Darcy} with 
the classical interface conditions \eqref{eq:IC-mass}, 
\eqref{eq:IC-momentum}, \eqref{eq:IC-BJJ} to the location of the sharp 
fluid-porous interface and to the Beavers--Joseph parameter $\alpha_\BJ$.

We evaluate the cross-sections at $x_1 = 0.4875$ (intersecting solid 
inclusions) and $x_1 = 0.5$ (between solid inclusions), where the flow is 
almost parallel to the interface and near the boundary at 
$x_1 = 0.025$, where the flow is almost perpendicular. The locations of 
the sharp fluid-porous interface for the macroscale model are $x_2 = 
-0.01$, $x_2 = 0$ (no shift) and  $x_2 = 0.01$.

\begin{figure}[!ht]
\begin{minipage}[t]{0.47\textwidth}
\centering
\includegraphics[width=1.0\textwidth]{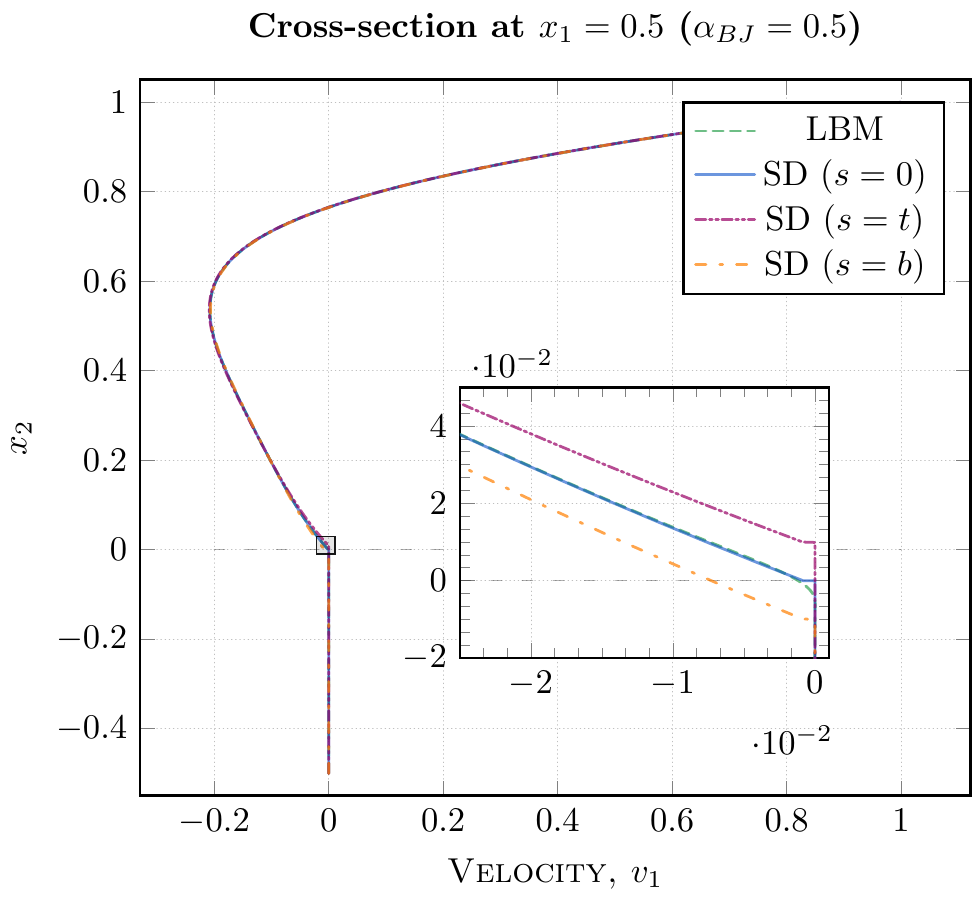} 
\caption{Velocity profiles ($x_1$-component) for the lid driven case at 
$x_1=0.5$ for different interface locations zoomed near the interface.}
\label{fig:LD-1}
\end{minipage}%
\hspace{0.725cm}
\begin{minipage}[t]{0.47\textwidth}
\centering
\includegraphics[width=1.0\textwidth]{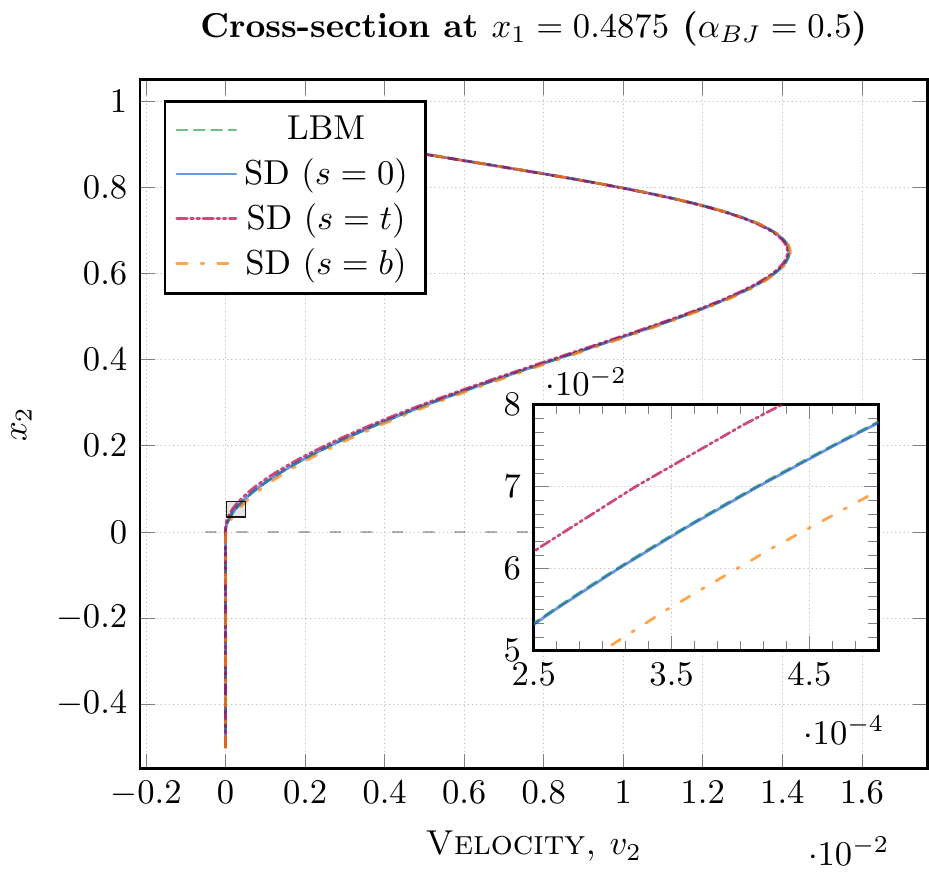} 
\caption{Velocity profiles ($x_2$-component) for the lid driven 
case at $x_1=0.4875$ for different 
interface locations zoomed near the interface.}
\label{fig:LD-2}
\end{minipage}
\end{figure}

Velocity profiles ($x_1$- and $x_2$-components) in and close to the 
horizontal centre of the domain are presented in Figures~\ref{fig:LD-1}, 
\ref{fig:LD-2}. The Beavers--Joseph parameter is 
taken $\alpha_\BJ = 0.5$ and 
different interface locations are considered: $(0,1) \times \{-0.01\}$ 
(shift below, $s=b$), $(0,1) \times \{0\}$ 
(no shift, $s=0$) and $(0,1) \times \{0.01\}$ 
(shift above, $s=t$). The profiles differ in the near-interface region and 
the best fit of the macroscale model to the pore-scale resolved 
simulations is obtained for $s=0$. This means that the sharp fluid-porous 
interface $\Gamma$ should be placed directly on top of the first row of 
solid inclusions (Fig.~\ref{fig:lid-driven}).

\begin{figure}[!h]
\centering
\includegraphics[width=0.5\textwidth]{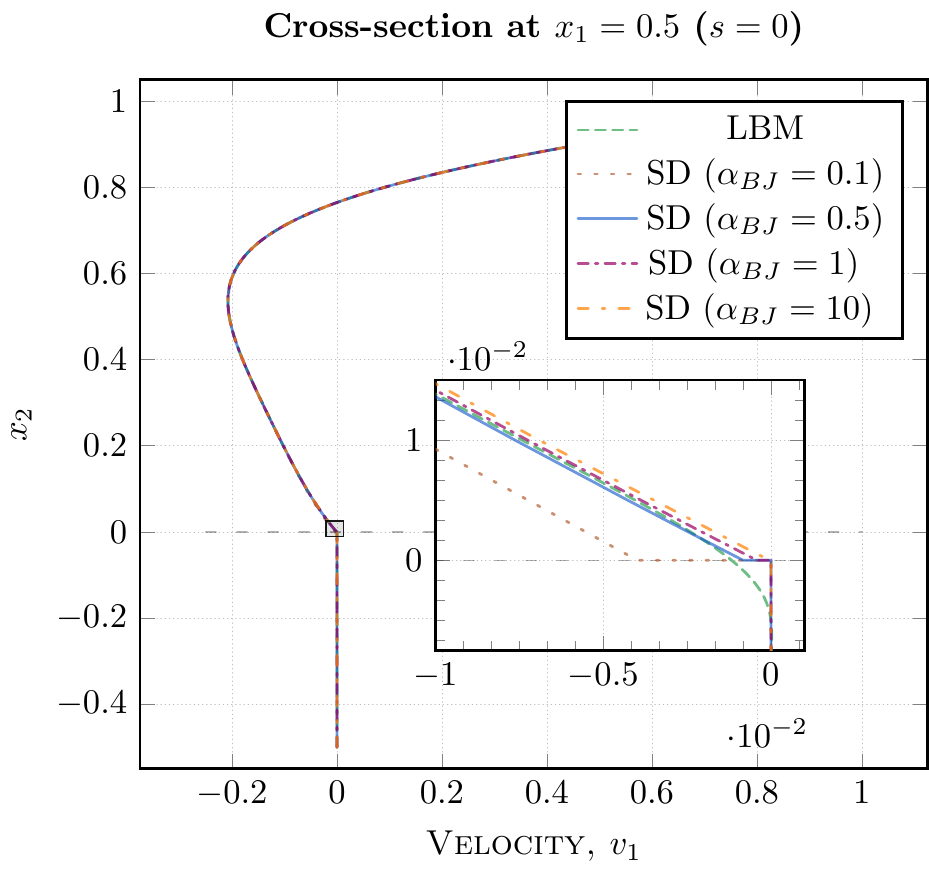} 
\caption{Velocity profiles ($x_1$-component) for the lid driven case at 
$x_1=0.5$ for different values of the 
Beavers--Joseph coefficient zoomed near the interface.}
\label{fig:LD-3}
\end{figure}

Figure~\ref{fig:LD-3} provides velocity profiles for different values of 
the Beavers--Joseph parameter $\alpha_\BJ$. The best fit is obtained for 
$\alpha_\BJ = 0.5$. This will be also confirmed later using interface 
conditions \eqref{eq:IC-Mik-parallel} obtained by homogenisation 
(Fig.~\ref{fig:LD-6}--\ref{fig:LD-a}). In 
the literature, $\alpha_\BJ = 1$ is usually used, which is not the optimal 
choice.

\begin{figure}[!ht]
\begin{minipage}[t]{0.47\textwidth}
\centering
\includegraphics[width=1.0\textwidth]{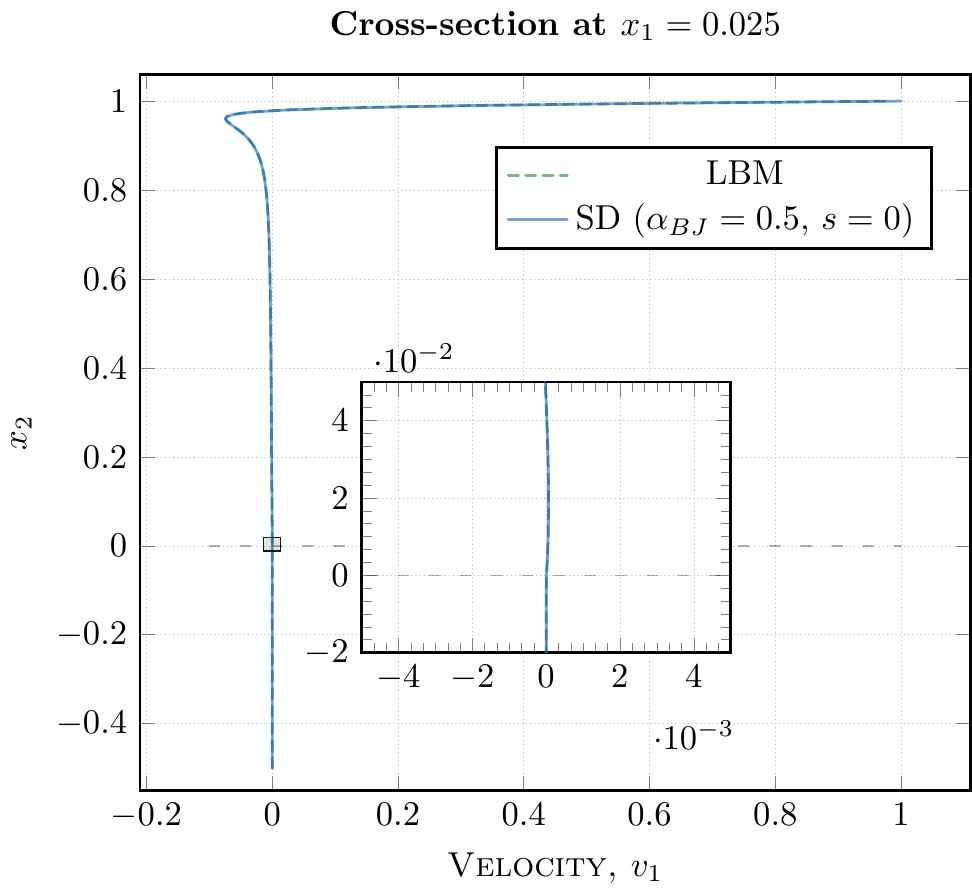} 
\caption{Velocity profiles ($x_1$-component) for the lid driven case at 
$x_1=0.025$ zoomed near the interface.}
\label{fig:LD-4}
\end{minipage}
\hspace{0.725cm}
\begin{minipage}[t]{0.47\textwidth}
\centering
\includegraphics[width=1.0\textwidth]{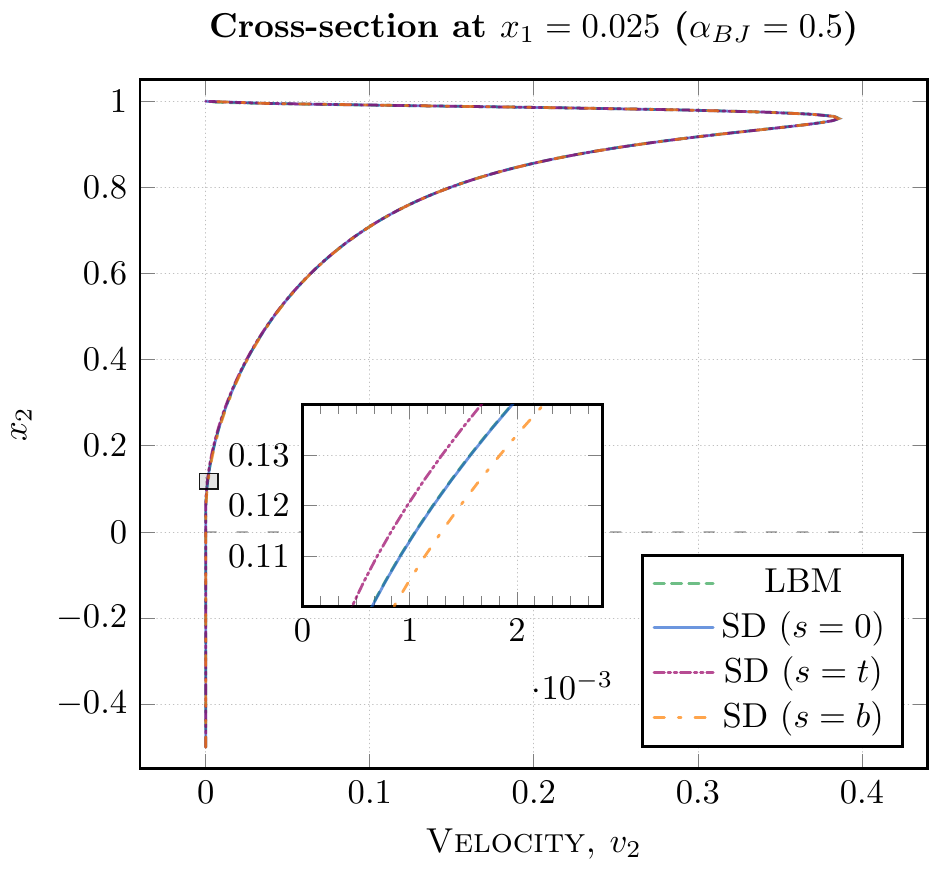} 
\caption{Velocity profiles ($x_2$-component) for the lid driven case at $x_1=0.025$ for different 
interface locations zoomed near the interface.}
\label{fig:LD-5}
\end{minipage}
\end{figure}

Velocity profiles near the left boundary are provided in 
Figures~\ref{fig:LD-4}, \ref{fig:LD-5}. Here, the location of the sharp 
interface and also the value of the Beavers--Joseph parameter $\alpha_\BJ$ 
has almost no effect on the $x_1$-component of the velocity 
(Fig.~\ref{fig:LD-4}). This is due to the fact that the velocity is 
almost perpendicular to the porous layer and the magnitude of the 
horizontal velocity is very small. However, for 
the vertical component of the velocity, the sharp interface location is 
important (Fig.~\ref{fig:LD-5}).

\subsubsection*{(ii) Interface conditions based on 
homogenisation}\label{sec:val-hom}

In this section, we study the sensitivity of the coupled macroscale 
Stokes--Darcy model \eqref{eq:1p1c-FF-mass}, \eqref{eq:1p1c-FF-momentum}, 
\eqref{eq:PM-mass}, \eqref{eq:PM-Darcy} with the interface conditions 
\eqref{eq:IC-Mik-parallel} proposed in \cite{Jaeger_etal_01, 
Jaeger_Mikelic_09}. 

\begin{figure}[!ht]
\begin{minipage}[t]{0.47\textwidth}
\centering
\includegraphics[width=1\textwidth]{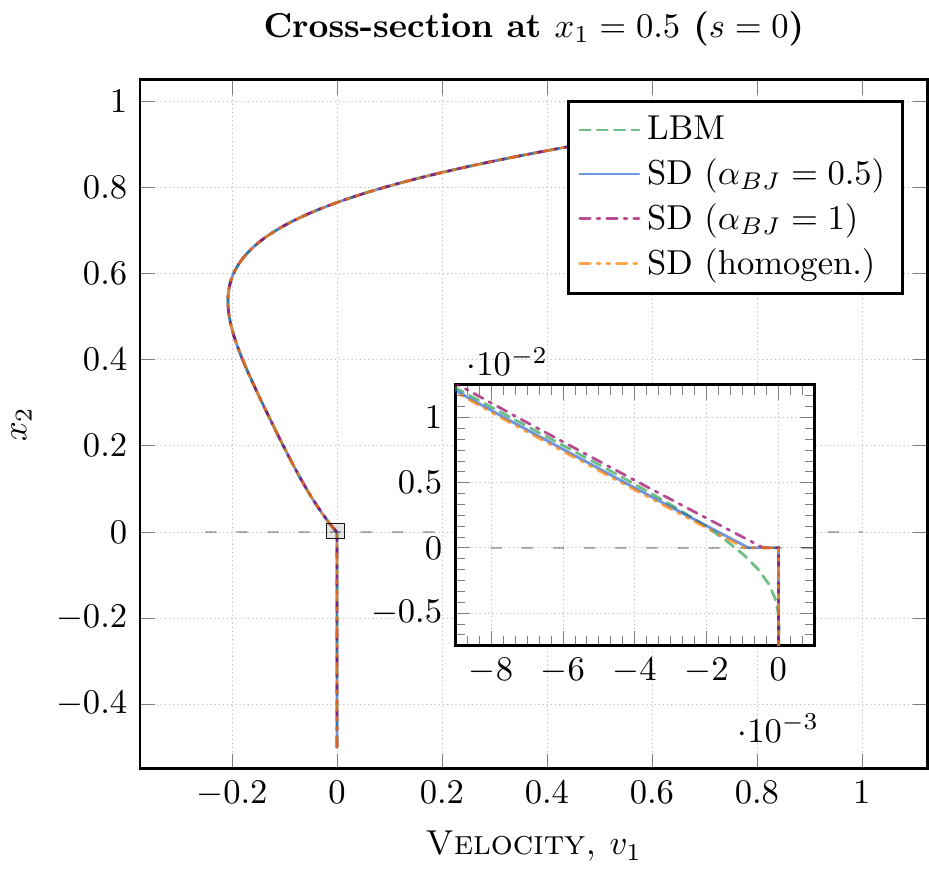} 
\caption{Velocity profiles ($x_1$-component) for the lid driven case at 
$x_1=0.5$ for classical interface conditions and conditions derived via 
homogenisation.}
\label{fig:LD-6}
\end{minipage}
\hspace{0.725cm}
\begin{minipage}[t]{0.47\textwidth}
\centering
\includegraphics[width=1\textwidth]{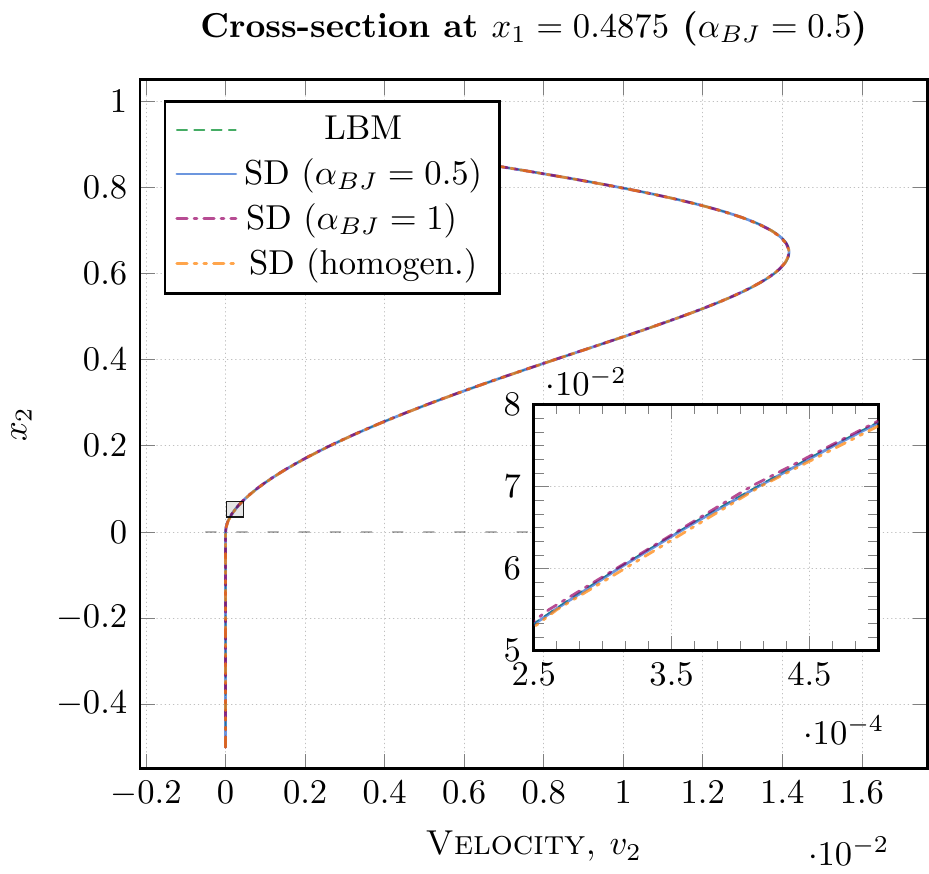} 
\caption{Velocity profiles ($x_2$-component) for the lid driven case at 
$x_1=0.4875$ for classical interface conditions and conditions derived via 
homogenisation.}
\label{fig:LD-7}
\end{minipage}
\end{figure}

To apply these interface conditions we need to 
compute the boundary layer constants for our geometry 
(Fig.~\ref{fig:lid-driven}). For these constants, we obtain the following values
\begin{equation}\label{eq:lbm-BL-constans}
C_1^{\text{bl}} = -5.41 \cdot 10^{-2}, \quad  C_\omega^{\text{bl}} = 0.
\end{equation} 

In Figures~\ref{fig:LD-6}, \ref{fig:LD-7} we provide velocity profiles for 
the lid driven cavity problem in and close to the horizontal centre of the 
domain. For the 
macroscale model we observed that the interface conditions derived by 
homogenisation (profile: SD (homogen.)) fit the best to the classical 
coupling conditions for $\alpha_\BJ = 0.5$. For such flow problems, we 
advise to use condition \eqref{eq:IC-Mik-parallel} instead of trying to 
find the optimal value for $\alpha_\BJ$.

Since the Beavers--Joseph--Saffman condition is the interface condition 
for the horizontal component of the velocity, the vertical component is 
only poorly sensitive to a change in the value of the Beavers--Joseph 
parameter $\alpha_\BJ$. We observed that the interface condition for the 
$x_2$-component of the velocity, given by the 
first equation in \eqref{eq:IC-Mik-parallel}, is nevertheless reasonable 
for such flow problems.

\begin{figure}[!h]
\centering
\includegraphics[width=0.47\textwidth]{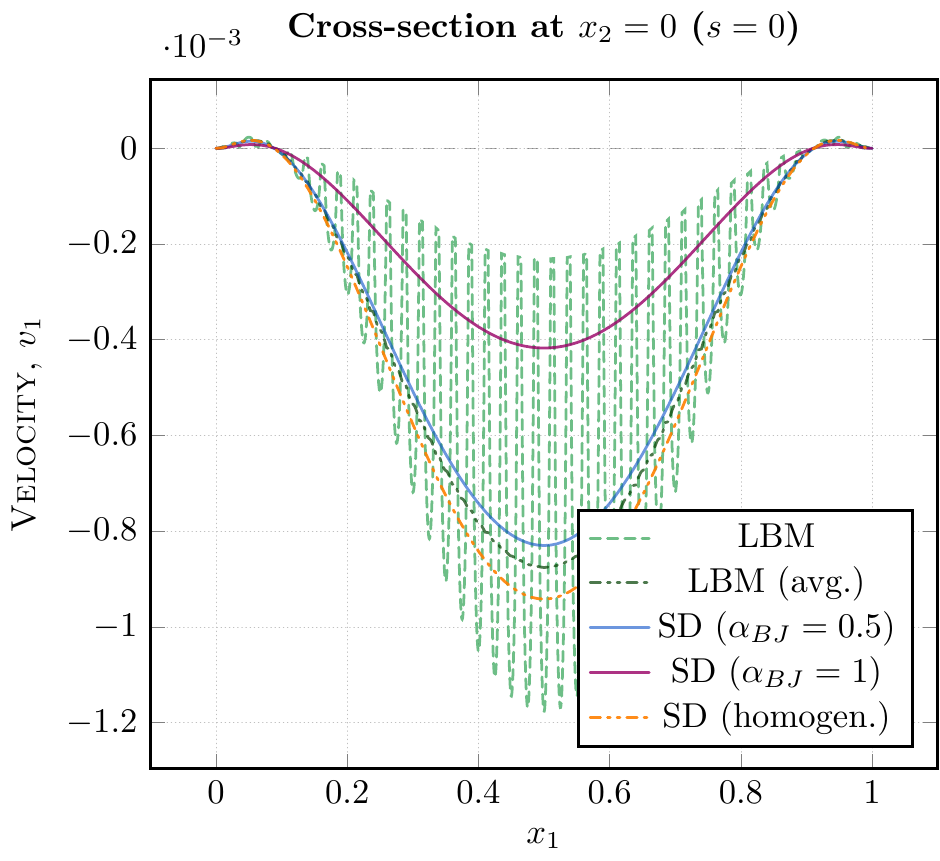} 
\caption{Velocity profiles ($x_1$-component) for the lid driven case at 
the sharp fluid--porous interface ($x_2=0$) for different values of the 
Beavers--Joseph coefficient $\alpha_\BJ$ and for conditions derived via 
homogenisation.}
\label{fig:LD-a}
\end{figure}

To demonstrate the sensitivity of the coupled Stokes--Darcy model to the 
choice of the Beavers--Joseph parameter $\alpha_\BJ$ in \eqref{eq:IC-BJJ}, 
we make cross-sections of the horizontal component of the velocity at the 
fluid-porous interface ($x_2 = 0$). The numerical simulation results for 
different values of $\alpha_\BJ$ and also for the interface conditions 
derived by homogenisation \eqref{eq:IC-Mik-parallel} are presented in 
Figure~\ref{fig:LD-a}.

For the macroscale simulations, the free-flow 
velocity at the interface is plotted.
The averaged pore-scale simulations (Fig.~\ref{fig:LD-a}, profile: LBM (avg.)) suggest that the choice $\alpha_\BJ=1$ which is commonly used in the literature is not 
optimal one. Again, as shown in Figures~\ref{fig:LD-3}, ~\ref{fig:LD-6}, 
~\ref{fig:LD-7}, we observed that $\alpha_\BJ = 0.5$ as well as the 
interface conditions \eqref{eq:IC-Mik-parallel} provide a much better 
fit.

\subsubsection{Staggered geometrical configuration}\label{sec:lid-2}

\begin{figure}[!h]
\centering
\includegraphics[width=0.5\textwidth]{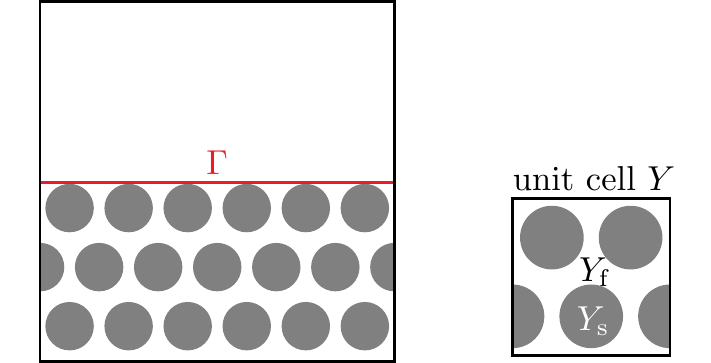}
\caption{Schematic illustration of the microscale setting and the 
location of the sharp fluid-porous interface $\Gamma$ (left) and the 
unit cell (right) for the staggered lid driven cavity and the infiltration problem.}
\label{fig:infiltration}
\end{figure}

Since the Beavers--Joseph parameter $\alpha_\BJ$ contains the information 
about the geometry of the interfacial region, we consider another 
geometrical configuration for the porous bed. 
In this case, the porous medium consists of solid inclusions that are 
arranged in a periodic but staggered manner (Fig.~\ref{fig:infiltration}).
The porosity is kept the same as for the channelised setting $\phi = 
0.4$. The number of solid inclusions is also the same: 40 inclusions in 
$x_1$-direction and 20 
inclusions in $x_2$-direction.

For this porous-medium geometry, we obtain the effective permeability 
tensor $\ten K = \varepsilon^2 \ten{\widetilde{K}}$ in the same way as before with
\begin{equation}\label{eq:k2}
\ten{\widetilde{K}} =
\begin{pmatrix}
1.325\cdot 10^{-4} & 0 \\
0 & 1.325\cdot 10^{-4}
\end{pmatrix}, 
\end{equation}
where $\varepsilon = 1 /20$ and the unit cell $Y$ contains several solid 
inclusions (Fig.~\ref{fig:infiltration}, right).

We analysed the macroscale problem \eqref{eq:1p1c-FF-mass}, 
\eqref{eq:1p1c-FF-momentum}, \eqref{eq:PM-mass}, \eqref{eq:PM-Darcy} with 
the classical interface conditions \eqref{eq:IC-mass}, 
\eqref{eq:IC-momentum}, \eqref{eq:IC-BJJ} and the appropriate boundary 
conditions at the external boundary.

\begin{figure}[!h]
\centering
\includegraphics[width=0.47\textwidth]{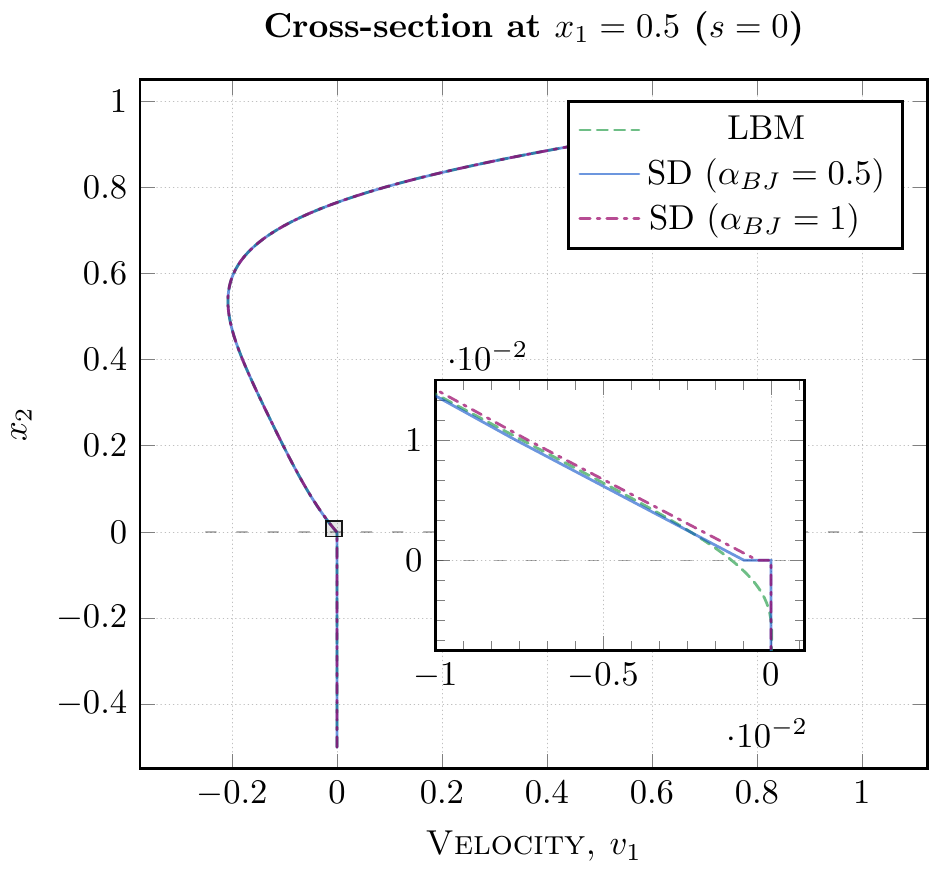} 
\caption{Velocity profiles ($x_1$-component) for the staggered lid driven case 
at $x_1=0.5$ for classcial interface conditions zoomed near the interface.}
\label{fig:lid-stagg}
\end{figure}

Figure~\ref{fig:lid-stagg} provides velocity profiles of the horizontal 
component of the velocity in the middle of the domain. Again, as we 
observed 
for the channelised geometry, the macroscale simulation results with 
$\alpha_\BJ=0.5$ fit much better to the pore-scale resolved simulations 
than the results with $\alpha_\BJ=1$ that is commonly used in the 
literature. 

\subsection{Infiltration}\label{sec:infiltration}

For the infiltration problem we consider only the staggered geometrical configuration of the porous medium (Fig.~\ref{fig:infiltration}) to 
avoid a channel flow in the vertical direction between columns of solid 
inclusions (Fig.~\ref{fig:lid-driven}). Therefore, the permeability tensor 
is identical with the one given in \eqref{eq:k2}.

The boundary conditions for the macroscale model given by equations 
\eqref{eq:1p1c-FF-mass}, \eqref{eq:1p1c-FF-momentum}, \eqref{eq:PM-mass}, 
\eqref{eq:PM-Darcy}, \eqref{eq:IC-mass}, \eqref{eq:IC-momentum}, 
\eqref{eq:IC-BJJ} for the infiltration problem read 
\begin{equation}\label{eq:infilt-bc}
\begin{split}
&\vec v_{\FF} = (0,-0.1 \sin (\pi x_{1})) \hspace*{0.35mm}  \textnormal{on 
} \{x_2 = 1\},\\
&\vec v_{\FF} = (0,0) \hspace*{19.25mm}  \textnormal{on } (\{0 \} \cup 
\{1\} ) \times (0,1), \\
&p_{\PM} = 100 \hspace*{20.15mm}  \textnormal{on } \{x_2 = -0.5\},\\
&\vec v_{\PM} \vdot \vec n_\PM = 0 \hspace*{15.85mm}  \textnormal{on } 
(\{0 
\} \cup \{1\} ) \! \times \!(-0.5,0).
\end{split}
\end{equation}

Here, we define the Reynolds number
\begin{equation}
Re=\frac{v_{\text{in,max}}\,\mathcal{L}}{\nu},
\end{equation}
where $v_{\text{in,max}}$ is the absolute maximum value of the inflow 
velocity ($x_2$-component). Similar as for the lid driven cavity problem 
(Sect.~\ref{sec:lid-driven}), in the LBM simulation, the whole domain is 
resolved by $2400 \times 3600$ 
lattice cells. The solid inclusions have a diameter of around 52 lattice cells. To resemble \eqref{eq:infilt-bc}, we 
use 
simple bounce back boundary conditions at the domain-boundaries in $x_1$-direction, at the inflow and at the solid 
inclusions. At the bottom of the domain in $x_2$-direction, we impose pressure boundary conditions based on the 
anti-bounce-back approach~\cite{Krueger_LBM} and set the lattice density to $\rho=1$. Again, $Re=0.01$ is chosen in 
order to 
approximate Stokes flow.

We make cross-sections for the velocity and pressure in the horizontal 
centre of the domain at $x_1=0.5$ and consider different locations of the 
sharp fluid-porous interface $\Gamma$ for the macroscale model. As in 
Section~\ref{sec:lid-driven}, the sharp interface is located at $x_2 = 
-0.01$ ($s=b$), 
$x_2 = 0$ (no shift, $s=0$) and $x_2 = 0.01$ ($s=t$).

\begin{figure}[!ht]
\begin{minipage}[t]{0.47\textwidth}
\centering
\includegraphics[width=1\textwidth]{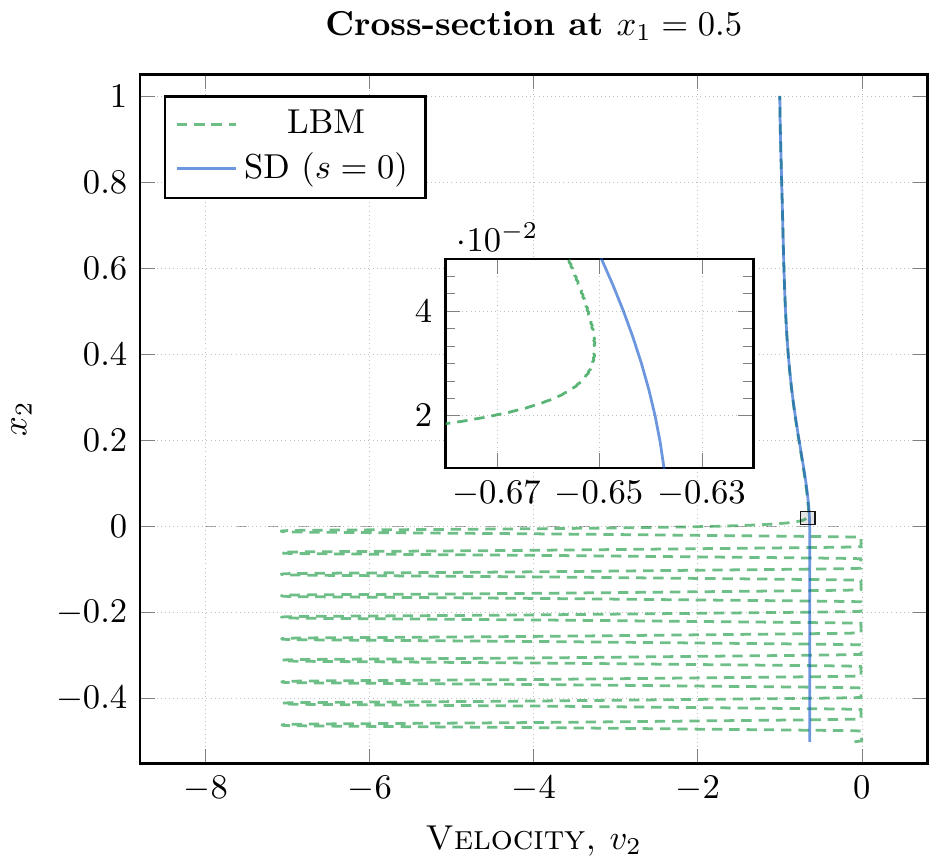} 
\caption{Velocity profiles ($x_2$-component) for the infiltration problem 
at $x_1=0.5$ zoomed near the interface.}
\label{fig:Inf-1}
\end{minipage}
\hspace{0.725cm}
\begin{minipage}[t]{0.47\textwidth}
\centering
\includegraphics[width=1\textwidth]{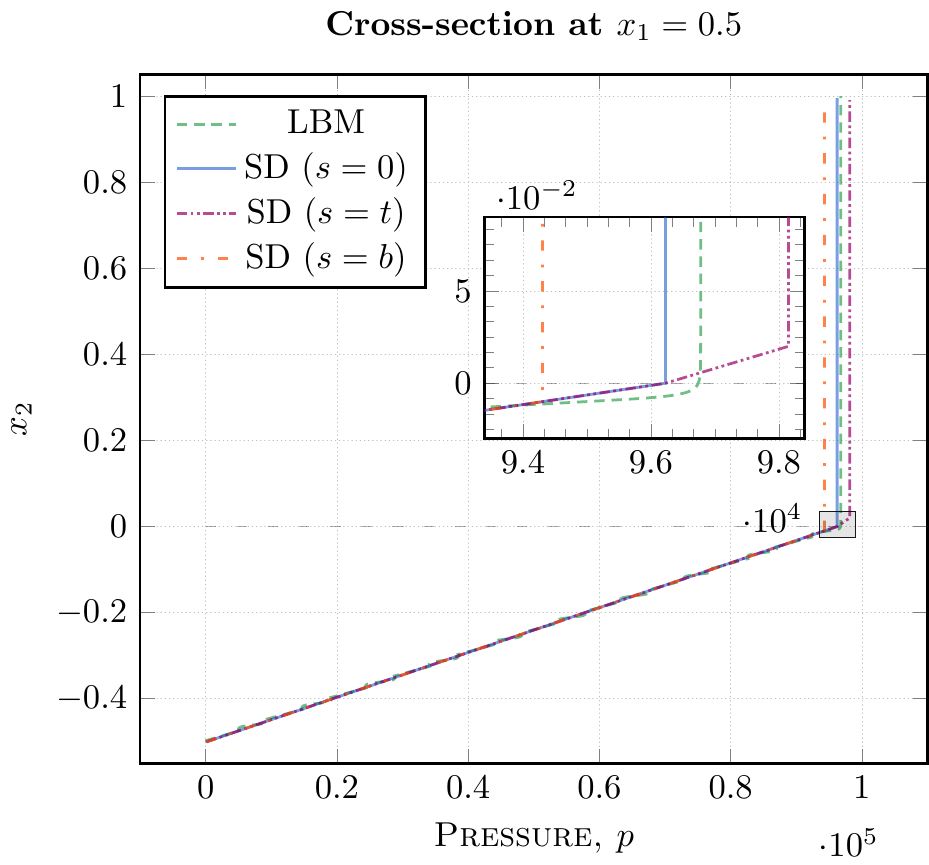} 
\caption{Pressure profiles for the infiltration problem at $x_1=0.5$ for 
different interface locations zoomed near the interface.}
\label{fig:Inf-2}
\end{minipage}
\end{figure}

In Figure~\ref{fig:Inf-1} we present the vertical component of the 
velocity computed for the macroscale model (SD) and the LBM simulations. 
We 
normalised the velocity by its maximal absolute 
value $v_{\text{in,max}}=0.1$. Although 
we used for the macroscale simulations the interface conditions 
\eqref{eq:IC-mass}, \eqref{eq:IC-momentum}, 
\eqref{eq:IC-BJJ}, which are valid for parallel flows, we get a quite good 
fit to the results of the pore-scale resolved models. This is due to the 
fact that the $x_1$-component of the velocity is almost zero. 
Therefore, we do not provide the profiles for the horizontal velocity 
component.

For the infiltration problem, we observed that the pressure field is very 
sensitive to the location of the sharp interface (Fig.~\ref{fig:Inf-2}).
Again, we obtain the best fit of the macroscale to the pore-scale 
model (Fig.~\ref{fig:Inf-2}, $s=0$) for the 
interface location presented in Figure~\ref{fig:infiltration}.
However, the differences between the pressure field obtained from the LBM 
simulations and the one from the best fitting interface location of the 
macroscopic model are more significant compared to the differences we 
observed for the velocity profiles.
Here, it is important to recall that the permeability tensor for the 
macroscale model is computed using the finite element framework 
\textsc{FreeFEM++}. We suspect that one cause of the more notable 
differences in the pressure profiles could be slight differences in the 
representation of the solid inclusions in the computational grids of LBM 
and \textsc{FreeFEM++}.

\section{Conclusion}\label{sec:conclusion}

In this paper, we validated two sets of interface conditions for 
the coupled Stokes--Darcy problem: (i) the classical conservation of 
mass across the interface, the balance of normal forces and the 
Beavers--Joseph--Saff\-man interface condition, and (ii) the interface 
conditions derived by means of the homogenisation and boundary 
layer theory~\cite{Jaeger_Mikelic_09, Jaeger_etal_01}. We considered 
different geometrical configurations (channelised and staggered 
arrangement of solid grains) and different flow problems (lid 
driven cavity over a porous bed and infiltration problem). All effective 
model parameters are computed numerically for the geometrical settings 
examined in the paper. To be able to compute these parameters 
using homogenisation theory, we considered periodic porous media. The 
interface conditions are validated and calibrated numerically by comparing 
macroscale simulation results against pore-scale simulations. As a 
pore-scale resolved 
model we used large scale lattice Boltzmann simulations.

Numerical simulation results demonstrate sensitivity of the coupled 
Stokes--Darcy problem to the choice of the effective model parameters, the 
location of the sharp interface and the interface conditions. We have 
tested a number of characteristic flow regimes in low Reynolds number cases 
and in all of them $\alpha_\BJ=0.5$ is proved to be the best candidate for 
the Beavers--Joseph parameter. We observed that the commonly 
used parameter $\alpha_\BJ = 1$ in the Beavers--Joseph--Saffman interface 
condition \eqref{eq:IC-BJJ} is often not the optimal choice. Moreover, 
the sharp fluid-porous interface should be located directly on the top of 
the first row of solid inclusions.

\begin{acknowledgements}
The authors thank Ivan Yotov, Jim Magiera and Christoph Rettinger for the 
valuable discussions. 
\newline
I. Rybak and E. Eggenweiler thank the Deutsche Forschungsgemeinschaft 
(DFG, German Research 
Foundation) for supporting this work by funding SFB 1313, Project Number 
327154368.\newline
C. Schwarzmeier and U. R\"ude thank the DFG for supporting project RU 422/27 and the Bundesministerium f\"ur Bildung 
und Forschung (BMBF, Federal Ministry of Education and 
Research) for supporting project SKAMPY (01IH15003A). The authors are grateful to the Regionales Rechenzentrum Erlangen 
(www.rrze.de) for providing access to supercomputing facilities.
\end{acknowledgements}

\bibliographystyle{spmpsci}      
\bibliography{LBM}  
\end{document}